%

\magnification=\magstep1
\def\forces{\parallel\!\!\! -}
\def\restrict{{\restriction}}
\def\Veskip{\vskip0.8truecm}
\def\Smallskip{\vskip1.4truecm}
\def\Bigskip{\vskip2.2truecm}
\def\Hoskip{\hskip1.6truecm}

\def\qed{{\vcenter{\hrule height.4pt \hbox{\vrule width.4pt height5pt
 \kern5pt \vrule width.4pt} \hrule height.4pt}}}
\def\ok{\vbox{\hrule height 8pt width 8pt depth -7.4pt
    \hbox{\vrule width 0.6pt height 7.4pt \kern 7.4pt \vrule width 0.6pt height 7.4pt}
    \hrule height 0.6pt width 8pt}}
\def\nt{{\leq}\kern-1.5pt \vrule height 6.5pt width.8pt depth-0.5pt \kern 1pt}
\def\sd{{\times}\kern-2pt \vrule height 5pt width.6pt depth0pt \kern1pt}
\def\notin{{\in}\kern-5.5pt / \kern1pt}
\def\ZZ{{\rm Z}\kern-3.8pt {\rm Z} \kern2pt}
\def\NN{{\rm I\kern-1.6pt {\rm N}}}

\def\DD{{\Bbb D}}

\def\BB{{\Bbb B}}

\def\PP{{\Bbb P}}
\def\RR{{\Bbb R}}
\def\KK{{\rm I\kern-1.6pt {\rm K}}}
\def\11{{\rm 1}\kern-2.2pt {\rm \vrule height6.1pt
    width.3pt depth0pt} \kern5.5pt}
\def\zp#1{{\hochss Y}\kern-3pt$_{#1}$\kern-1pt}

\def\egs{\vrule height 6pt width.5pt depth 2.5pt \kern1pt}
\font\small=cmr8 scaled\magstep0
\font\smalli=cmti8 scaled\magstep0
\font\capit=cmcsc10 scaled\magstep0
\font\capitg=cmcsc10 scaled\magstep1

\font\dunhg=cmr10 scaled\magstep1
\font\dunhgg=cmr10 scaled\magstep2

\font\sanse=cmss10 scaled\magstep0
\overfullrule=0pt
\openup1.5\jot
\input mssymb

\centerline{}
\Bigskip
\centerline{\dunhgg Combinatorial properties of Hechler forcing}
\footnote{}{{\openup-6pt {\small {\smalli 1991 
Mathematics subject classification.}
Primary 03E40, Secondary 03E15 28A05 54H05 \par {\smalli Key words and
phrases.} Hechler real, category, measure, random real, Cicho\'n's
diagram, almost disjoint family\endgraf}}}
\Bigskip
\centerline{\capitg J\"org Brendle$^{1,}$\footnote{$^*$}{{\small
The first author would like to thank the MINERVA-foundation
for supporting him}}, Haim Judah$^{1,}$\footnote{$^{**}$}{{\openup-7pt
\small
The second and third authors would like to thank the Basic Research
Foundation (the Israel Academy of Sciences and Humanities) for
partially supporting them\endgraf}} and Saharon Shelah$^{2,**}$}
\Smallskip
{\baselineskip=0pt {\small
\noindent $^1$
Abraham Fraenkel Center for Mathematical Logic,
Department of Mathematics,
Bar--Ilan University,
52900 Ramat--Gan, Israel 
\smallskip
\noindent $^2$ Institute of Mathematics,
The Hebrew University,
Jerusalem, Israel}}
\Bigskip
\centerline{\capit Abstract}
\bigskip
\noindent Using a notion of rank for Hechler forcing we show:
1) assuming $\omega_1^V = \omega_1^L$, there is no real 
in $V[d]$ which is eventually different from the reals in
$L[d]$, where $d$ is Hechler over $V$; 2) adding one Hechler real
makes the invariants on the left-hand side of Cicho\'n's
diagram equal $\omega_1$ and those on the right-hand side 
equal $2^\omega$ and produces a maximal almost
disjoint family of subsets of $\omega$ of size
$\omega_1$; 3) there is no perfect set of random
reals over $V$ in $V[r][d]$, where $r$ is random over $V$ and
$d$ Hechler over $V[r]$, thus answering a question of the
first and second authors.
\vfill\eject
\noindent{\dunhg Introduction}
\Smallskip
In this work we use a notion of rank first introduced by James
Baumgartner and Peter Dordal in [BD, $\S$ 2] and later developed
independently by the third author in [GS, $\S$ 4] to show that
{\it adding a Hechler real} has strong combinatorial consequences.
Recall that the {\it Hechler
p. o. $\DD$} is defined as follows.
\smallskip
\centerline{$(s,f) \in \DD \Longleftrightarrow s \in \omega^{< \omega}
\; \land \; f \in \omega^\omega \; \land \; s \subseteq f \; 
\land \; f$ strictly increasing}
\smallskip
\centerline{$(s,f)  \leq (t,g) \Longleftrightarrow s \supseteq t \;
\land \;
\forall n \in \omega \; (f(n) \geq g(n))$}
\smallskip
\noindent We note here that our definition differs from the usual one
in that it generically adds a strictly increasing function from
$\omega$ to $\omega$. This is, however, a minor point making the definition
of the rank in section 1 easier. We indicate at the end of $\S$ 1
how it can be changed to get the corresponding results in $\S\S$ 2 and
4 for {\it classical} Hechler forcing.
\par
The theorems of section 2 are all consequences of one technical result
which is expounded in 2.1. We shall sketch how some changes in the latter's
argument prove that adding one Hechler real produces a maximal almost
disjoint family of subsets of $\omega$ of size $\omega_1$ (2.2.).
Recall that $A, B \subseteq \omega$ are said to be {\it almost
disjoint} ({\it a. d.} for short) iff $\vert A \cap B \vert <
\omega$; ${\cal A} \subseteq [\omega]^\omega$ is an {\it a. d.
family} iff the members of ${\cal A}$ are pairwise a. d.; and
${\cal A}$ is a {\it m. a. d. family} ({\it maximal almost
disjoint family}) iff it is a. d. and maximal with this property. ---
We shall then show that assuming $\omega_1^V =
\omega_1^L$, there is no real in $V[d]$ which is eventually different
from the reals in $L[d]$, where $d$ is Hechler over $V$ (2.4.).
Here, we say that given models $M \subseteq N$ of $ZFC$, a real
$f \in \omega^\omega \cap N$ is {\it eventually different} from the
reals in $M$ iff $\forall g \in \omega^\omega \cap M \; \forall^\infty
n \; (g(n) \neq f(n))$, where $\forall^\infty n$ abbreviates {\it
for all but finitely many $n$}. (Similarly, $\exists^\infty n$ will
stand for {\it there are infinitely any $n$}.) ---
Next we will prove that adding one Hechler real makes the invariants
on the left-hand side of Cicho\'n's diagram equal $\omega_1$
and those on the right-hand side equal $2^\omega$ (2.5.). These
invariants (which describe combinatorial properties of measure
and category on the real line, and of the eventually dominating
order on $\omega^\omega$) will be defined, and the shape of
Cicho\'n's diagram explained, in the discussion preceding
the result in $\S$ 2. Theorem 2.5. should be seen as a continuation of
research started by Cicho\'n and Pawlikowski in [CP] and [Pa].
They investigated the effect of {\it adding a Cohen or a random real}
on the invariants in Cicho\'n's diagram. ---
We close section 2 with an application concerning absoluteness
in the projective hierarchy (2.6.); namely we show that {\it $\Sigma_4^1 -
\DD$-absoluteness} (which means that $V$ and $V[d]$, where $d$ is 
Hechler over $V$, satisfy the same $\Sigma_4^1$-sentences with
parameters in $V$) implies that $\omega_1^V > \omega_1^{L[r]}$ for
any real $r$; in particular $\omega_1^V$ is inaccessible in $L$.
So, for projective statements, Hechler forcing is much stronger 
than Cohen or random forcing for $\Sigma_n^1$-Cohen-absoluteness
($\Sigma_n^1$-random-absoluteness) is true in any model gotten
by adding $\omega_1$ Cohen (random) reals [Ju, $\S$ 2].
\par
In $\S$ 3 we leave Hechler forcing for a while to deal with
{\it perfect sets of random reals} instead, and to continue
a discussion initiated in [BaJ] and [BrJ]. Recall that given
two models $M \subseteq N$ of $ZFC$, we say that $g \in 
\omega^\omega \cap N$ is a {\it dominating real over
$M$} iff $\forall f \in \omega^\omega \cap M \;
\forall^\infty n \; (g(n) > f(n))$;
and $r \in 2^\omega \cap N$ is {\it random over $M$}
iff $r$ avoids all Borel null sets coded in $M$ iff $r$
is the real determined by some filter which is $\BB$-generic
over $M$ (where $\BB$ is the algebra of Borel sets
of $2^\omega$ modulo the null sets ({\it random
algebra}) -- see [Je, section 42] for details).
--- A tree $T \subseteq 2^{< \omega}$ is {\it perfect} iff $\forall t 
\in T \; \exists s \supseteq t \; (s \hat{\;} \langle 0
\rangle \in T \; \land s \hat{\;} \langle 1 \rangle \in T)$.
For a perfect tree $T$ we let $[T] := \{ f \in 2^\omega ;
\; \forall n \; (f \restrict n \in T ) \}$ denote the
set of its branches. Then $[T]$ is a perfect set (in the
topology of $2^\omega$). Conversely, given a perfect
set $S \subseteq 2^\omega$ there is perfect tree $T \subseteq
2^{< \omega}$ such that $[T] = S$. This allows us to confuse
perfect sets and perfect trees in the sequel; in particular,
we shall use the symbol $T$ for both the tree and the set of its branches.
--- We will show in 3.1. that given models $M \subseteq N$
of $ZFC$ such that there is a perfect set of random reals
in $N$ over $M$, either there is a dominating real in $N$
over $M$ or $\mu(2^\omega \cap M) = 0$ in $N$. This result is sharp and
has some consequences concerning the relationship between
cardinals related to measure and to the eventually dominating
order on $\omega^\omega$ (cf [BrJ, 1.9] and the discussion
preceding 3.2. for details).
\par
The argument for theorem 3.1. together with the techniques
of $\S$ 1 yield the main result of section 4; namely,
there is no perfect set of random reals over $M$ in
$M[r][d]$, where $r$ is random over $M$, and $d$ Hechler
over $M[r]$ (4.2.). This answers questions 2 and 2' in [BrJ].
\bigskip
{\sanse Notation.} Our notation is fairly standard. We refer
the reader to [Je] and [Ku] for set theory in general
and forcing in particular. 
\par
Given a finite sequence $s$ (i.e. either
$s \in 2^{<\omega}$ or $s \in \omega^{<\omega}$), we let $lh(s) 
:= dom (s) $ denote the length of $s$; for $\ell \in lh(s)$,
$s \restrict \ell$ is the restriction of $s$ to $\ell$. 
$\hat{\;}$ is used for concatenation of sequences; and
$\langle \rangle$ is the empty sequence. 
Given a perfect tree $T \subseteq 2^{< \omega}$ and $s \in T$,
we let $T_s := \{ t \in T ; \; t \subseteq s$ or $s \subseteq t
\}$. --- 
Given a p.o. $\PP \in V$, we
shall denote $\PP$-names by symbols like $\tau$, 
$\breve f$, $\breve T$, ...
and their interpretation in $V[G]$ (where $G$ is $\PP$-generic
over $V$) by $\tau [G]$, $\breve f [G]$, $\breve T [G]$, ... 
\bigskip
{\sanse Acknowledgement.} 
We would like to thank Andrzej Ros{\l}anowski for several helpful
discussions.
\Bigskip

\noindent{\dunhg $\S$ 1. Prelude --- a notion of rank for Hechler
forcing}
\Smallskip
{\sanse 1.1.} {\capit Main Definition} (Shelah, see [GS, $\S$ 4]
--- cf also [BD, $\S$ 2]).
Given $t \in \omega^{< \omega}$ strictly increasing 
and $A \subseteq \omega^{<
\omega}$, we define by induction when the {\it rank}
$rk (t, A)$ is $\alpha$. \par
\item{(a)} $rk(t,A) =0$ iff $t \in A$. \par
\item{(b)} $rk (t,A) = \alpha$ iff  for no $\beta < \alpha$
we have $rk (t,A) = \beta$, but there are $m \in \omega$ and
$\langle t_k ; \; k \in \omega \rangle$ such that $\forall
k \in \omega$: $t \subseteq t_k$, $t_k \in \omega^m$,
$t_k (lh(t)) \geq k$, and $rk(t_k , A) < \alpha$. $\qed$ \par
\smallskip
Clearly, the rank is either $< \omega_1$ or undefined
(in which case we say $rk = \infty$). We repeat the proof of
the following result for it is the main tool for $\S\S$ 2
and 4.
\smallskip
{\sanse 1.2.} {\capit Main Lemma} (Baumgartner--Dordal [BD, $\S$ 2]
and Shelah [GS, $\S$ 4]).
{\it Let $I \subseteq \DD$ be dense. Set $A: = \{ t ; \;
\exists f \in \omega^\omega$ such that $(t,f) \in I \}$.
Then $rk(t^*,A) < \omega_1$ for any $t^* \in \omega^{< \omega}$.}
\smallskip
{\it Proof.} Suppose $rk(t^*, A) = \infty$ for some $t^* \in
\omega^{< \omega}$. Let $S: = \{ s \in \omega^{<\omega}$ strictly
increasing$; \; t^* \subseteq s$ and for all $s^*$ with
$t^* \subseteq s^*$ and with $\forall i \in dom(s^*) \setminus dom(t^*)
\; (s^* (i) \geq s(i) )$, we have $rk(s^* ,A) = \infty \}$.
$S \subseteq \omega^{< \omega}$ is a tree with stem $t^*$.
\par
Suppose $S$ has an infinite branch $\langle s_i ; \; i \in \omega 
\rangle$ (i.e. $s_0 = t^*$, $lh(s_i) = lh(t^*) + i$, and $s_i 
\subseteq s_{i+1}$). Let $g$ be the function defined by this
branch: $g = \bigcup_{i \in \omega} s_i$. Then $(t^* , g)
\in \DD$. Choose $(t,f) \leq (t^*, g)$ such that $(t,f)
\in I$. Then $t \in A$, i.e. $rk(t,A) = 0$; but also $t \in S$,
i.e. $rk(t,A) = \infty$, a contradiction. 
\par
So suppose $S$ has no infinite branches, and let $s^*$ be a
maximal point in $S$. Then we have a sequence $\langle t_k ;
\; k \in \omega \rangle$ such that $lh (t_k) = lh (s^*) + 1$,
$t_k(lh(s^*)) \geq k$, $t^* \subseteq t_k$, $\forall i \in
dom(s^*) \setminus dom (t^*) \; (t_k (i) \geq s^* (i))$,
and $rk(t_k ,A) < \infty$. Now we can find a subset $B \subseteq
\omega$ and $lh (t^*) \leq m \leq lh(s^*)$ and $t \in \omega^m$
such that $\forall k \in B \; (t_k \restrict m = t)$ and
$k < \ell$, $k, \ell \in B$, implies $t_k (lh(t)) < t_\ell (lh (t))$.
Hence the sequence $\langle t_k ; \; k \in B \rangle$ witnesses
$rk(t,A) < \infty$. On the other hand $t \in S$; i.e. $rk(t,A)
= \infty$, again a contradiction. $\qed$
\bigskip
Usually Hechler forcing $\DD '$ is defined as follows.
\smallskip
\centerline{$(s,f) \in \DD ' \Longleftrightarrow s \in \omega^{<\omega}
\; \land \; f \in \omega^\omega \; \land \; s \subseteq f$}
\smallskip
\centerline{$(s,f) \leq (t,g) \Longleftrightarrow s \supseteq t \; \land
\; \forall n \in \omega \; (f(n) \geq g(n))$}
\smallskip
\noindent We sketch how to introduce a rank on $\DD '$ having the
same consequences as the one on $\DD$ defined above. Let $\Omega
= \{ t ; \; dom(t) \subseteq \omega \; \land \; \vert t \vert < \omega
\; \land \; rng(t) \subseteq \omega \}$. Given $t \in \Omega$ and
$A \subseteq \omega^{<\omega}$ we define by induction when the {\it rank} $rk
(t,A)$ is $\alpha$.
\par
\item{(a)} $rk(t,A) = 0$ iff $t \in A$. \par
\item{(b)} $rk(t,A) = \alpha$ iff for no $\beta < \alpha$ we have
$rk(t,A) = \beta$, but there are $M \in [\omega]^{<\omega}$ and
$\langle t_k ; \; k \in \omega \rangle$ such that $dom(t)
\subset M$ and $\forall k \in \omega$: $ t \subseteq t_k$,
$t_k \in \omega^M$, $rk(t_k , A) < \alpha$ and $\forall i \in
M \setminus dom(t) \; \forall k_1 \neq k_2 \; (t_{k_1} (i) \neq
t_{k_2} (i) )$. \par
\noindent We leave it to the reader to verify that the result corresponding to
1.2. is true for this rank on $\DD '$, and that the theorems of
$\S\S$ 2 and 4 can be proved for $\DD '$ in the same way as they are
proved for $\DD$.
\Bigskip

\noindent{\dunhg $\S$ 2. Application I --- the effect of adding one 
Hechler real on the invariants in Cicho\'n's diagram}
\Smallskip
Before being able to state the main result of this section
(the consequences of which will be 1) and 2) in the abstract)
we have to set up some notation.
\par
Let ${\cal A} \subseteq [\omega]^\omega$
be an a. d. family. 
We will produce a set of $\DD$-names
$\{ \tau_A ; \; A \in {\cal A} \}$ for functions in $\omega^\omega$
as follows. For each $A \in {\cal A}$ fix $f_A : A \to
\omega$ onto with $\forall n \;\exists^\infty m\in A \;
(f_A (m) = n)$. Now, if $r \in \omega^\omega$ is a real
having the property that $\{ n \in \omega ; \; r(n) \in A \}$
is infinite, let $g_r : \omega \to \omega$ be an enumeration
of this set (i.e. $g_r (0) := $ the least $n$ such that $r(n) \in A$;
$g_r (1) :=$ the least $n > g_r (0)$ such that $r(n) \in A$; etc.).
In this case we let $\tau_A (r) : \omega \to \omega$ be defined as
follows.
\smallskip
\centerline{$\tau_A(r) (n) := f_A (r (g_r (n)))$.}
\smallskip
\noindent As $A$ is infinite, we have $\forces_\DD " \vert rng(\breve d)
\cap A \vert = \omega "$, where $\breve d$ is the name for the
Hechler real; in particular $\tau_A (d)$ will be defined in the 
generic extension. Thus we can think of $\langle \tau_A (\breve d) ;
\; A \in {\cal A} \rangle$ as a sequence of names in Hechler forcing
for objects in $\omega^\omega$.
\bigskip
{\sanse 2.1.} {\capit Main Theorem.} {\it Whenever ${\cal A} \subseteq
[\omega]^\omega$ is an a. d. family in the ground model $V$, $d$ is
Hechler over $V$, and $f \in \omega^\omega$ is any real in $V[d]$,
then $\{ A \in {\cal A} ; \; \forall^\infty n \; (f(n) \neq \tau_A
(d) (n)) \}$ is at most countable (in $V[d]$).}
\smallskip
{\it Remark.} Slight changes in the proof show that, in fact,
$\{ \tau_A ; \; A \in {\cal A} \}$ is a Luzin set in $V[d]$
for uncountable ${\cal A}$. (Recall that an uncountable set of reals is
called {\it Luzin} iff for all meager sets $M$, $M \cap S$ is
at most countable.)
\smallskip
{\it Proof.} The proof uses the main lemma (1.2.) as principal tool.
Let $\breve f$ be a $\DD$-name for a real (for an element of
$\omega^\omega$). Let $I_n$ be the set of conditions deciding
$\breve f \restrict (n+1)$ ($n \in \omega$). All $I_n$ are dense.
Let $D_n := \{ t ; \; \exists f \in \omega^\omega $ such that $(t,f)
\in I_n \}$ (cf the main lemma). We want to define when a
set $A \in {\cal A}$ is {\it $n$-bad}. \par
For each $t \in \omega^{< \omega} \setminus
D_n$ strictly increasing we can find
(according to the main lemma for $D_n$) an $m \in \omega$ and
$\langle t_k ; \; k \in \omega \rangle$ such that for all $k \in
\omega$: $t_k$ is strictly increasing, $t \subseteq t_k$, $t_k
\in \omega^m$, $t_k(lh(t)) \geq k$, and $rk(t_k, D_n) < rk(t,D_n)$.
Let $m_t := m - lh(t)$. We define by induction on $i < m_t$ when
$A \in {\cal A}$ is {\it $t-i-n$-bad}. Along the way we also
construct sets $B_i$ ($i < m_t$). \par
$i = 0$. Let $B_0 = \omega$. If there is $A \in {\cal A}$ 
such that $A \cap \{ t_k (lh(t)) ; \; k \in B_0 \}$ is infinite,
choose such an $A_0$ and let $A_0$ be $t-0-n$-bad. Now let $B_1 =
\{ k \in \omega ; \; t_k(lh(t)) \in A_0 \}$. If there is no such
$A$, let $B_1 = B_0 = \omega$. \par
$i \to i+1$ ($i+1 < m_t$). We assume that $B_{i+1}$ is defined and
infinite. If there is  $A \in {\cal A}$ such that $A \cap \{
t_k(lh(t) +i+1) ; \; k \in B_{i+1} \}$ is infinite, choose
such an $A_{i+1}$ and let $A_{i+1}$ be $t-(i+1)-n$-bad. Now
let $B_{i+2} = \{ k \in B_{i+1} ; \; t_k (lh(t)+i+1) \in A_{i+1} \}$.
If there is no such $A$, let $B_{i+2} = B_{i+1}$. \par
In the end, we set $B_t := B_{m_t}$. We say that $A \in {\cal A}$ is
{\it $n$-bad} iff it is $t-i-n$-bad for some strictly increasing
$t \in \omega^{<\omega} \setminus
D_n$ and $i < m_t$. Finally $A \in {\cal A}$ is
{\it bad} iff it is $n$-bad for some $n \in \omega$. Let ${\cal A}_{\breve
f} = \{ A \in {\cal A} ; \; A$ bad $\}$. Since for $n \in
\omega$, $t \in \omega^{<\omega}$ and $i < m_t$ at most
one $A \in {\cal A}$ is $t - i -n$-bad, ${\cal A}_{\breve f}$
is countable.
\smallskip
{\capit Claim.} {\it If $A \in {\cal A} \setminus {\cal A}_{\breve f}$,
then $\forces_\DD \exists^\infty n  \; ( \breve f (n) = \tau_A (\breve d)
(n) )$.}
\smallskip
{\it Remark.} Clearly this claim finishes the proof of the main theorem.
\smallskip
{\it Proof.} Suppose not, and choose $(s,g) \in \DD$, $k \in \omega$,
and $A \in {\cal A} \setminus {\cal A}_{\breve f}$ such that
\smallskip
\centerline{$(s,g) \forces_\DD \forall n \geq k \; (\breve f (n) \neq
\tau_A (\breve d) (n))$.}
\smallskip
\noindent Let $\ell \geq k$ be such that $\vert rng(s) \cap A \vert \leq \ell$;
i.e. $s$ {\it does not decide} the value of $\tau_A(\breve d) (\ell)$.
By increasing $s$, if necessary, we can assume that $\vert rng(s) \cap
A \vert = \ell$. Let $Y:= \{ t \in \omega^{< \omega} ; \; t $
strictly increasing, $s \subseteq t$, $\forall i \in dom(t) \setminus
dom(s) \; (t(i)  \geq g(i))$, and $
\vert rng(t) \cap A \vert = \ell \}$. Choose $t \in Y$ such that
$rk(t,D_\ell)$ is minimal.
\smallskip
{\capit Subclaim.} {\it $rk(t,D_\ell) = 0$.}
\smallskip
{\it Proof.} Suppose not. Then choose by the main lemma (1.2.)
$m \in \omega$ and $\langle t_k ; \; k \in \omega \rangle$ (i.e.
all $t_k$ are strictly increasing, $t \subseteq t_k$,
$t_k \in \omega^m$, $t_k(lh(t)) \geq k$, and $rk(t_k ,D_\ell) < rk(t,
D_\ell)$).
In fact, we require that $m$ and $\langle t_k ; \; k \in \omega 
\rangle$ are the same as the ones chosen for $\ell$, $t$ in the 
definition of $\ell$-badness.
Let $m_t = m - lh(t)$ as above, and look at $B_t$. By construction
(as $A$ is not $t-i-\ell$-bad for any $i < m_t$) and almost-disjointness,
$A \cap \{ t_k (lh(t) + i) ; \; k \in B_{i + 1} \}$ is finite for
all $i < m_t$. So there is $k \in B_t$ such that $rng(t_k) \cap A
= rng (t) \cap A$, i.e. $\vert rng(t_k) \cap A \vert = \ell$, and
$t_k(i) \geq g(i)$ for all $i \in dom(t_k) \setminus dom(s)$. Hence
$t_k \in Y$ and $rk(t_k, D_\ell) < rk(t,D_\ell)$, contradicting
the minimality of $rk(t,D_\ell)$. This proves the subclaim. $\qed$
\smallskip
{\it Continuation of the proof of the claim.} As $rk(t,D_\ell) = 0$
we have an $h \in \omega^\omega$ such that $(t,h) \in I_\ell$.
Then $(t,max(h,g)) \leq (s,g)$, and this condition decides the value 
of $\breve f$ at $\ell$ without deciding the value of $\tau_A (\breve d)$
at $\ell$. Suppose that $(t,max(h,g)) \forces_\DD " \breve f (\ell)
= j " $. Now choose $i \geq max(h,g) (lh(t))$ such that $i \in A$
and $f_A(i) =j$ (this exists by the choice of the function $f_A$).
Then
\smallskip
\centerline{$(t \hat{\;} \langle i \rangle, max(h,g)) \forces_\DD
\breve f (\ell) = j = f_A (i) = f_A (\breve d (g_{\breve d} (\ell))) =
\tau_A (\breve d) (\ell)$.}
\smallskip
\noindent This final contradiction ends the proof of the claim and
of the main theorem. $\qed$ $\qed$
\bigskip
We will sketch how a modification of this argument gives the
following result.
\smallskip
{\sanse 2.2.} {\capit Theorem.} {\it After adding one Hechler real
$d$ to $V$, there is a maximal almost disjoint family of subsets
of $\omega$ of size $\omega_1$ in $V[d]$.}
\smallskip
{\it Sketch of proof.} We start with an observation which will
relate Luzin sets and maximal almost disjoint families.
\smallskip
{\capit Observation.} {\it Let $\langle N_\alpha ; \; \omega
\leq \alpha < \omega_1 \rangle$, $\langle h_\alpha ; \; \omega
\leq \alpha < \omega_1 \rangle$ and $\langle r_\alpha ; \;
\omega \leq \alpha < \omega_1 \rangle$ be sequences such
that $N_\alpha \prec H(\kappa)$ is countable and $N_\alpha
\prec N_\beta$ for $\alpha < \beta$, $h_\alpha \in \alpha^\omega
\cap N_\alpha$ is one-to-one and onto, $r_\alpha \in \omega^\omega$
is Cohen over $N_\alpha$ and $\langle r_\alpha ; \; \alpha <
\beta \rangle \in N_\beta$. Define recursively sets $C_\alpha$
for $\alpha < \omega_1$. $\langle C_n ; \; n \in \omega \rangle
$ is a partition of $\omega$ into countable pieces lying in $N_\omega$.
For $\alpha \geq\omega$, $C_\alpha := \{ r_\alpha (n) ; \; n \in
\omega \; \land \; \forall m < n \; (r_\alpha (n) \not\in C_{h_\alpha
(m)} ) \}$. Then $\{ C_\alpha ; \; \alpha \in \omega_1 \}$
is an a. d. family.}
\smallskip
{\it Proof.} The construction gives almost--disjointness. So it suffices
to show that each $C_\alpha$ is infinite. But this follows from the fact 
that each $r_\alpha$ is Cohen over $N_\alpha$ and that the union of
finitely many $C_\beta$'s (for $\beta < \alpha$) is coinfinite. $\qed$
\smallskip
Now let ${\cal A} = \langle A_\alpha ; \; \alpha < \omega_1 \rangle \in
V$ be an a. d. family. As $\langle \tau_{A_\alpha} (d) ; \; \alpha
< \omega_1 \rangle$ is Luzin in $V[d]$ (see the remark following the
statement of theorem 2.1.) we can find a strictly increasing function 
$\phi : \omega_1 \setminus \omega \to \omega_1$ and sequences 
$\langle N_\alpha ; \; \omega \leq \alpha < \omega_1 \rangle$,
$\langle h_\alpha ; \; \omega \leq \alpha < \omega_1 \rangle$
such that for $r_\alpha := \tau_{A_{\phi(\alpha)}} (d)$ the requirements
of the above observation are satisfied. By $ccc$--ness of $\DD$, we
may assume that $\phi \in V$; and hence, that $\phi = id$, thinning ${\cal
A}$ out if necessary. We want to show that the resulting family
$\langle C_\alpha ; \; \alpha < \omega_1 \rangle$ is a m. a. d.
family. \par
For suppose not. Then there is a $\DD$-name $\breve C$ such that
\smallskip
\centerline{$\forces_\DD \forall \alpha < \omega_1 \; ( \vert
\breve C_\alpha \cap \breve C \vert < \omega)$.}
\smallskip
\noindent Let $\breve f$ be the $\DD$-name for the strictly increasing
enumeration of $\breve C$. As in the proof of 2.1. we let $I_n$ be
the set of conditions deciding $\breve f \restrict (n+1)$, $D_n :=
\{ t ; \; \exists f \in \omega^\omega \; ((t,f) \in I_n ) \}$, and
we define when a set $A \in {\cal A}$ is $n$-bad (so that at most
countably many sets will be $n$-bad). \par
Furthermore, for each $\alpha < \omega_1$ we let $\sigma_\alpha$ be the
$\DD$-name for a natural number such that 
\smallskip
\centerline{$\forces_\DD \breve C_\alpha \cap \breve C \subseteq
\sigma_\alpha$.}
\smallskip
\noindent We let $I_\alpha '$ be the set of conditions deciding
$\sigma_\alpha
$, $D_\alpha ' := \{ t ; \; \exists f \in \omega^\omega \; ((t,f) \in
I_\alpha ' ) \}$; analogously to the proof of theorem 2.1. we define when
a set $A \in {\cal A}$ is $\alpha$-bad (so that at most countably
many sets will be $\alpha$-bad). \par
Next choose $\alpha < \omega_1$ such that \par
1) if $A_\beta$ is $n$-bad for some $n$, then $\beta < \alpha$; \par
2) if $\beta < \alpha$ and $A_\gamma$ is $\beta$-bad, then
$\gamma < \alpha$. \par
\smallskip
{\capit Claim.} {\it $\forces_\DD \vert \breve C_\alpha \cap \breve C \vert
= \omega$.}
\smallskip
{\it Proof.} Suppose not, and choose $(s,g) \in \DD$ and $k \in \omega$ 
such that \smallskip
\centerline{$(s,g) \forces_\DD \breve C_\alpha \cap \breve C \subseteq k$.}
\smallskip
\noindent Let $\ell \geq k$ be such that $\vert rng(s) \cap A_\alpha \vert
\leq \ell$; without loss $\vert rng(s) \cap A_\alpha \vert = \ell$.
Let $Y := \{ t \in \omega^{< \omega} ; \; t$ strictly increasing, $s
\subseteq
t$, $\forall i \in dom(t) \setminus dom(s) \; (t(i) \geq g(i)) $, and
$\vert rng(t) \cap A_\alpha \vert = \ell \}$. By the argument of the
subclaim in the proof of 2.1. there is a $t \in Y$ such that $\forall
m < \ell \; (rk(t,D_{h_\alpha (m)} ' ) = 0 )$. Hence there is an
$h \in \omega^\omega$ such that $(t,h) \in \bigcap_{m < \ell} I_{h_\alpha
(m)} '$. Without loss $h \geq g$. Then $(t,h) \leq (s,g)$, and
this condition decides the values of $\sigma_{h_\alpha (m)}$ ($m<\ell$);
suppose that $(t,h) \forces_\DD " \forall m < \ell \; (\sigma_{h_\alpha (m)}
=s_m)"$. Choose $\ell '$ larger that the maximum of the $s_m$ ($m<\ell$)
and $k$. Again using the argument of the subclaim (2.1.) find $t' 
\supseteq t$ such that $\forall i \in dom(t') \setminus dom (t) \;
(t'(i) \geq h(i))$, $\vert rng(t') \cap A_\alpha \vert = \ell$, and
$rk(t', D_{\ell '} ) = 0$. Thus there exists an $h' \in \omega^\omega$
such that
\smallskip
\centerline{$(t',h') \forces_\DD " \breve f (\ell ') = j "$ for
some $j$.}
\smallskip
\noindent Without loss $h' \geq h$. Then $(t',h') \leq (t,h)$. As 
$\forces_\DD " \breve f$ is strictly increasing", $j \geq \ell ' \geq k$;
by construction we have in particular that $(t',h') \forces_\DD
"\forall m < \ell \; (j \not\in \breve C_{h_\alpha (m)})"$. Choose
$i \geq h' (lh (t'))$ such that $i \in A_\alpha$ and $f_{A_\alpha} (i)
= j$. Then \smallskip
\centerline{$(t' \hat{\;} \langle i \rangle , h') \forces_\DD  \breve f
(\ell ') = j = f_{A_\alpha} (i) = \tau_{A_\alpha} (\breve d) (\ell) =
\breve r_\alpha (\ell) \in \breve C_\alpha$.}
\smallskip
\noindent This final contradiction proves the claim, and the theorem
as well. $\qed$ $\qed$
\smallskip
In our proof we constructed a m. a. d. family of size $\omega_1$ from
a Luzin set in $V[d]$. We do not know whether this can be done in $ZFC$.
\smallskip
{\sanse 2.3.} {\capit Question} (Fleissner, see [Mi, 4.7.]) {\it
Does the existence of a Luzin set imply the existence of a m. a. d. family
of size $\omega_1$?}
\smallskip
{\it Remark.} It is consistent that there is a m. a. d. family of
size $\omega_1$, but no Luzin set. This is known to be true in
the model obtained by adding at least $\omega_2$ random reals to
a model of $ZFC + CH$.
\bigskip
We next turn to consequences of theorem 2.1.
\smallskip
{\sanse 2.4.} {\capit Theorem.} {\it Let $V \subseteq W$ be universes
of set theory, $\omega_1^V = \omega_1^W$. Then no real in $W[d]$
is eventually different from the reals in $V[d]$, where $d$ is Hechler
over $V$.}
\smallskip
{\it Remark.} Remember that Hechler forcing has an absolute definition.
So $d$ will be Hechler over $V$ as well.
\smallskip
{\it Proof.} Let ${\cal A} \subseteq [\omega]^\omega$ be an almost
disjoint family in $V$ of size $\omega_1$. Assume that the functions
$f_A$ for $A \in {\cal A}$ (defined at the beginning of this
section) are also in $V$. Then each real in $W[d]$ can
only be eventually different from countably many of the
reals in $\{ \tau_A (d) ; \; A \in {\cal A} \}
\in V[d]$, by the main theorem. $\qed$ 
\bigskip
To be able to explain our next corollary to the main theorem, we
need to introduce a few cardinals. Given a $\sigma$-ideal
${\cal I} \subseteq P(2^\omega)$, we let
\smallskip
\itemitem{$add({\cal I}) :=$} the least $\kappa$ such that 
$\exists {\cal F} \in [{\cal I}]^\kappa \; (\bigcup
{\cal F} \not\in {\cal I})$;
\par
\itemitem{$cov ({\cal I}) :=$} the least $\kappa$ such that
$\exists {\cal F} \in [{\cal I}]^\kappa \; (\bigcup
{\cal F} = 2^\omega)$;
\par
\itemitem{$unif({\cal I}) :=$} the least $\kappa$ such that
$[2^\omega]^\kappa \setminus {\cal I} \neq \emptyset$;
\par
\itemitem{$cof ({\cal I}) :=$} the least $\kappa$ such that
$\exists {\cal F} \in [{\cal I}]^\kappa \; \forall A \in {\cal I}
\; \exists B \in {\cal F} \; (A \subseteq B)$.
\smallskip
\noindent We also define
\smallskip
\itemitem{$b :=$} the least $\kappa$ such that $\exists {\cal F}
\in [\omega^\omega]^\kappa \; \forall f \in \omega^\omega
\; \exists g \in {\cal F} \; \exists^\infty n \;
(g(n) > f(n))$;
\par
\itemitem{$d :=$} the least $\kappa$ such that $\exists {\cal F}
\in [\omega^\omega]^\kappa \; \forall f \in \omega^\omega \;
\exists g \in {\cal F} \; \forall^\infty n \; (g(n) > f(n))$.
\smallskip
\noindent If ${\cal M}$ is the ideal of meager sets, and
${\cal N}$ is the ideal of null sets, then we can arrange
these cardinals in the following diagram (called Cicho\'n's
diagram).
\bigskip
\centerline{{\quad} \Hoskip $cov({\cal N})$ \Hoskip
$unif({\cal M})$ \Hoskip $cof({\cal M})$ \Hoskip $cof({\cal N})$
\Hoskip $2^\omega$}
\Veskip
\centerline{$b$ \hskip 2.5truecm $d$}
\Veskip
\centerline{$\omega_1$ \Hoskip $add({\cal N})$ \Hoskip
$add({\cal M})$ \Hoskip $cov({\cal M})$ \Hoskip $unif({\cal N})$
\Hoskip {\quad}}
\bigskip
\noindent (Here, the invariants grow larger, as one moves up and to the
right in the diagram.) The dotted line says that $add({\cal M})=
min\{ b, cov({\cal M}) \}$ and $cof({\cal M}) = max\{ d,
unif({\cal M}) \}$. For the results which determine the
shape of this diagram, we refer the reader to [Fr]. A survey
on independence proofs showing that no other relations can be
proved between these cardinals can be found in [BJS].
We shall need the following characterizations of the
cardinals $unif({\cal M})$ and $cov({\cal M})$, which are due to
Bartoszy\'nski [Ba].
\smallskip
$unif({\cal M}) =$ the least $\kappa$ such that $\exists {\cal F}
\in [\omega^\omega]^\kappa \; \forall g \in \omega^\omega \;
\exists f \in {\cal F} \; \exists^\infty n \; (f(n) = g(n))$;
\par
$cov({\cal M}) =$ the least $\kappa$ such that $\exists {\cal F}
\in [\omega^\omega]^\kappa \; \forall g \in \omega^\omega \;
\exists f \in {\cal F} \; \forall^\infty n \; (f(n) \neq g(n))$.
\smallskip
\noindent We are ready to give our next result, which says essentially
that after adding one Hechler real, the invariants on the left-hand
side of the above diagram 
all equal $\omega_1$, whereas those on the right-hand side
are all equal to $2^\omega$.
\smallskip
{\sanse 2.5.} {\capit Theorem.} {\it After adding one Hechler
real $d$ to $V$, $unif({\cal M}) = \omega_1$ and $cov({\cal M}) =
2^\omega$ in $V[d]$.}
\smallskip
{\it Proof.} (i) Let ${\cal A} \subseteq [\omega]^\omega$ be
an a. d. family of size $\omega_1$ in $V$. Then by the main
theorem no real is
eventually different from $\{ \tau_A (d) ; \; A \in {\cal A} \}$,
giving $unif({\cal M}) = \omega_1$ (by Bartoszy\'nski's characterization).
\par
(ii) Let ${\cal A} \subseteq [\omega]^\omega$ be an a. d.
family of size $2^\omega$ in $V$ (such a family exists, see
e.g. [Ku, chapter II, theorem 1.3]). Suppose $\kappa =
cov({\cal M}) < 2^\omega$, and let $\{ g_\alpha ; \; \alpha
< \kappa \}$ be a family of functions such that $\forall g \in
V[d] \cap \omega^\omega \; \exists \alpha < \kappa \; \forall^\infty
n \; (g(n) \neq g_\alpha (n))$, using Bartoszy\'nski's
characterization. As $\vert {\cal A} \vert =
2^\omega > \kappa$, there is ${\cal A}' \subseteq {\cal A}$,
$\vert {\cal A}' \vert \geq \omega_1$, and $\alpha < \kappa$ such
that $\forall A \in {\cal A}' \; \forall^\infty n \; (\tau_A
(d) (n) \neq g_\alpha (n))$. This contradicts the main theorem.
$\qed$
\smallskip
{\it Remark.} Instead of Bartoszy\'nski's characterization we
could have used the fact that $\{ \tau_A (d) ; \; A \in {\cal A} \}$ 
is a Luzin set (see the remark after 2.1.). We leave it to the reader 
to verify that the existence of a Luzin set implies $unif({\cal M})
= \omega_1$; and that the existence of a Luzin set of size $2^\omega$
implies $cov({\cal M}) = 2^\omega$.
\bigskip
We close with an application concerning absoluteness in the
projective hierarchy. We first recall a notion due to the
second author [Ju, $\S$ 2]. Given a universe of set theory $V$
and a forcing notion $\PP \in V$ we say that $V$ is {\it $\Sigma^1_n -
\PP$-absolute} iff for every $\Sigma^1_n$-sentence $\phi$ with
parameters in $V$ we have $V \models \phi$ iff $V^\PP \models \phi$.
So this is equivalent to saying that $\RR^V \prec_{\Sigma^1_n}
\RR^{V^\PP}$. Note that Shoenfield's Absoluteness Lemma [Je,
theorem 98] says that $V$ is alway $\Sigma^1_2 - \PP$-absolute.
Furthermore, $\Sigma^1_3 - \DD$- absoluteness is equivalent to
{\it all $\Sigma^1_2$-sets have the property of Baire} [Ju, $\S$
2]. This is a consequence of Solovay's classical characterization
of the latter statement which says that it is equivalent to:
{\it for all reals $a$, the set of reals Cohen over $L[a]$
is comeager.} 
\smallskip
{\sanse 2.6.} {\capit Theorem.} {\it $\Sigma_4^1-\DD$-absoluteness
implies that $\omega_1 > \omega_1^{L[r]}$ for any real $r$.}
\smallskip
{\it Proof.}
Suppose there is an $a \in \RR$ such that $\omega_1^{L[a]}
= \omega_1^V$. By $\Sigma^1_3 - \DD$-absoluteness we have that
all $\Sigma^1_2$-sets have the property of Baire (see above);
i.e. $\forall b \in \RR$ ($Co(L[b])$ is comeager) ($Co(M)$
denotes the set of reals Cohen over some model $M$ of $ZFC$).
Note that $x \in Co(L[b])$ is equivalent to
\smallskip
\centerline{$\forall c \; ( c \not\in L[b] \cap BC \; \lor
\; \hat c$ is not meager $\lor \; x \not\in \hat c)$,}
\smallskip
\noindent where $BC$ is the set of Borel codes which is $\Pi_1^1$
[Je, lemma 42.1], and for $c \in BC$, $\hat c$ is the set coded by
$c$. As $L[b]$ is $\Sigma_2^1$ [Je, lemma 41.1], $Co(L[b])$ is
a $\Pi^1_2$-set. Hence $\forall b \in \RR \; (Co(L[b])$ is
comeager$)$ which is equivalent to
\smallskip
\centerline{$\forall b \exists c \; ( c \in BC \; \land \;
\hat c$ is meager $\land \; \forall x \; (x \in \hat c
\; \lor \; x \in Co(L[b])))$}
\smallskip
\noindent is a $\Pi^1_4$-sentence. So it is true in $V^\DD$ by 
$\Sigma^1_4 - \DD$-absoluteness; in particular $Co(L[a][d])$
is comeager in $V[d]$ which implies that there is a dominating
real in $V[d]$ over $L[a][d]$, contradicting theorem 2.4.
$\qed$
\bigskip
{\sanse 2.7.} {\capit Question.} {\it Are there results similar
to theorems 2.4., 2.5., and 2.6. for Amoeba forcing or Amoeba-meager
forcing?}
\smallskip
We conjecture that the answer is yes because both the Amoeba
algebra and the Amoeba-meager algebra contain $\DD$ as a complete
subalgebra (see [Tr, $\S$ 6]; a definition of the algebras can
also be found there). But there doesn't seem to be a way
to introduce a rank on these algebras (as in $\S$ 1).
\Bigskip

\noindent{\dunhg $\S$ 3. Interlude --- perfect sets of random reals}
\Smallskip
{\sanse 3.1.} {\capit Theorem.} {\it Let $V \subseteq W$ be
models of $ZFC$. Suppose there is a perfect set of random
reals in $W$ over $V$. Then either \par
\item{1)} there is a dominating real in $W$ over $V$; or \par
\item{2)} $\mu(2^\omega \cap V) = 0$ in $W$.}
\smallskip
{\it Proof.} Suppose not, and let $T \in W$ be a perfect set
of random reals. Define $f \in \omega^\omega \cap W$ as follows.
\smallskip
\centerline{$f(i) = min \{ k ; \; \forall \sigma \in T \cap 2^i
\; (\vert T_\sigma \cap 2^k \vert > 4^i ) \}$}
\smallskip
\noindent Let $g \in \omega^\omega \cap V$ be such that 
$\exists^\infty i \; (g(i) \geq f(i))$. Let $U$ be the family of
all $u \in \prod_{i \in \omega} P(2^{g(i)})$ such that $u(i)
\subseteq 2^{g(i)}$ and ${\vert u(i) \vert \over 2^{g(i)}}
=2^{-i}$. $U$ can be thought of as a measure space
(namely, for $u \subseteq 2^{g(i)}$ with ${\vert u \vert \over
2^{g(i)}} = 2^{-i}$ let $\mu_i (u) = {1 \over {2^{g(i)} \choose
2^{g(i) - i} }}$; and let $\mu$ be the product measure of the
$\mu_i$). \par
Let $N \prec \langle H(\kappa)^W , ... \rangle$ be countable with
$g , T \in N$. As $\mu(2^\omega \cap V) \neq 0$ in $W$, we cannot
have that $2^\omega \cap V \subseteq \cup \{ B ; \;
\mu(B) = 0, \; B \in N, \; B$ Borel $ \}$; i.e. there are reals
in $V$ which are random over $N$. Let $u^* \in U$ be such a real.
Using $u^*$ we can define a measure zero set $B$ in $V$ as follows.
\smallskip
\centerline{$B = \{ h \in 2^\omega ; \; \exists^\infty i
\; (h \restrict g(i) \in u^* (i) ) \}$}
\smallskip
\noindent Let (for $k \in \omega$) $B_k = \{ h \in 2^\omega ;
\; \forall i \geq k \; (h \restrict g(i) \not\in u^* (i) ) \}$.
Clearly $2^\omega \setminus B = \cup_{k \in \omega} B_k$; and
the $B_k$ form an increasing chain of perfect sets of positive
measure. \par
As all reals in $T$ are random over $V$ we must have $T \subseteq
\cup_{k \in \omega} B_k$. This gives us $\sigma \in T$ and $k
\in \omega$ such that $T_\sigma \subseteq B_k$ (otherwise
choose $\sigma_0 \in T$ such that $\sigma_0 \not\in B_0$,
$\sigma_1 \in T_{\sigma_0}$ such that $\sigma_1 \not\in
B_1$, etc. This way we construct a branch in $T$ which does
not lie in $\cup_{k \in \omega} B_k$, a contradiction).
\par
By construction, we know that for infinitely many $i$, we have
$\vert T_\sigma \cap 2^{g(i)} \vert > 4^i$ and $u^* (i) 
\cap (T_\sigma \cap 2^{g(i)}) = \emptyset$. For each such $i$
and $u \subseteq 2^{g(i)}$ with ${\vert u \vert \over 2^{g(i)}}
= 2^{-i}$,
the probability that $u \cap (T_\sigma \cap 2^{g(i)})
= \emptyset$ (in the sense of the measure $\mu_i$ defined above) is
$$\leq ({2^{g(i)} - 4^i \over 2^{g(i)}})^{2^{g(i) - i}} \leq
(e^{-{4^i \over 2^{g(i)}}})^{2^{g(i) - i}} =
e^{-2^i}.$$
So the probability that this happens infinitely often is
zero. But $u^*$ is random over $N$, a contradiction. $\qed$
\smallskip
{\capit Corollary} (Cicho\'n [BaJ, $\S$ 2]). 
{\it If $r$ is random over $V$, then
there is no perfect set of random reals in $V[r]$ over $V$.} $\qed$
\smallskip
{\it Remark.} Theorem 3.1. is best possible in the following sense.
\par \item{1)} It is consistent that there are $V \subseteq W$ and a perfect
tree $T$ of random reals in $W$ over $V$ and $\mu^* (2^\omega \cap V) >
0$ in $W$ ($\mu^*$ denotes outer measure). To see this add a Laver real
$\ell$ to $V$ and then a random real $r$ to $V[\ell]$; set $W = V[\ell]
[r]$. By [BaJ, theorem 2.7] there is a perfect tree of random reals in
$W$ over $V$; and by [JS, $\S$ 1] $\mu^* (2^\omega \cap V) > 0$ in
$V[\ell]$ and hence in $W$. \par
\item{2)} It is consistent that there are $V \subseteq W$ and a perfect
tree $T$ of random reals in $W$ over $V$ and no dominating real in $W$
over $V$ (see [BrJ, theorem 1]).
\bigskip
Before being able to state some consequences of this result,
we need to introduce two further cardinals.
\smallskip
\itemitem{$wcov ({\cal N}) :=$} the least $\kappa$ such that
$\exists {\cal F} \in [{\cal N}]^\kappa \; (2^\omega \setminus\bigcup
{\cal F}$ does not contain a perfect set);
\par
\itemitem{$wunif ({\cal N}) :=$} the least $\kappa$ such that 
there is a family ${\cal F} \in [[2^{< \omega}]^\omega]^\kappa$
of perfect sets with $\forall N \in {\cal N} \; \exists
T \in {\cal F} \; (N \cap T = \emptyset)$.
\smallskip
\noindent We can arrange these cardinals and some of those
of the preceding section in the following diagram.
\bigskip
\centerline{$2^\omega$}
\Veskip
\centerline{$cof({\cal N})$}
\Veskip
\centerline{$cov({\cal N})$ \hskip 4truecm $d$ \hskip 4truecm $wunif({\cal N})$}
\Veskip
\centerline{$wcov({\cal N})$ \hskip 4truecm $b$ \hskip 4truecm $unif({\cal N})$}
\Veskip
\centerline{$add({\cal N})$}
\Veskip
\centerline{$\omega_1$}
\bigskip
\noindent (Here the invariants get larger as one moves up in the
diagram.) 
The dotted line says that
$wcov ({\cal N}) \geq min \{ cov({\cal N}) , b \}$ (and dually,
$wunif ({\cal N}) \leq max \{ unif ({\cal N}) , d \}$) (see
[BaJ, $\S$ 2] or [BrJ, 1.9]). Using the above result we get
\smallskip
{\sanse 3.2.} {\capit Theorem.} {\it (i) $wcov({\cal N}) \leq
max\{ b, unif ({\cal N}) \}$; \par
(ii) $wunif({\cal N}) \geq min \{ d , cov({\cal N}) \}$ --- In fact,
given $V \subseteq W$ models of $ZFC$ such that in $W$ there is a
real which is random over a real which is unbounded over $V$, 
there exists a null set $N \in W$ such that for all perfect sets
$T \in V$, $T \cap N \neq \emptyset$.}
\smallskip
{\it Proof.} (i) follows immediately from theorem 3.1; and the
first sentence of (ii) follows from the last sentence of (ii).
The latter is proved by an argument which closely follows the lines
of the proof of theorem 3.1, and is therefore left to the reader. $\qed$
\smallskip
The most interesting question concerning the relationship of
the cardinals in the above diagram is the following (question 3'
of [BrJ]).
\smallskip
{\sanse 3.3.} {\capit Question.} {\it Is it consistent that $wcov({\cal N})
> d$? Dually, is it consistent that $wunif({\cal N}) < b$?}
\Bigskip
\vfill\eject

\noindent{\dunhg $\S$ 4. Application II --- adding a Hechler real over
a random real does not produce a perfect set of random reals}
\Smallskip
{\sanse 4.1.} {\capit Theorem.} {\it Let $V \subseteq W$ be models 
of $ZFC$ such that \par
\item{1)} there is no dominating real in $W$ over $V$; \par
\item{2)} $2^\omega \cap V$ is non-measurable in $W$. \par
\noindent Then there is no perfect set of random reals in $W[d]$,
where $d$ is Hechler over $W$.}
\smallskip
{\it Remark.} This result clearly contains theorem 3.1. as a special
case; still we decided to bring the latter as a separate result
because it has consequences for the cardinals involved (see above, 3.2.).
Also, the proof of theorem 4.1. can be seen as a combination of the
argument for 3.1. and the techniques developed in $\S$ 1.
\smallskip
{\sanse 4.2.} {\capit Corollary.} {\it There is no perfect set of random
reals in $V [r] [d]$, where $r$ is random over $V$, and $d$ is
Hechler over $W=V[r]$.} $\qed$
\smallskip
{\it Proof of theorem 4.1.} We work in $W$. 
Let $\breve T$ be a $\DD$-name for a perfect 
tree. We want to show that $T = \breve T [G]$ ($G$ $\DD$-generic
over $W$) contains reals which are not random over $V$. 
We say that $A \subseteq \omega^{<\omega}$ is {\it large} iff 
$\forall (s,f) \in \DD \; \exists s' \in A $ with $(s' , f) \leq (s,f)$
(By $(s',f)$ we mean here and in the sequel the condition $(s',f')$
where $f' \restrict dom (s') = s' $ and $ f' (n) = f(n)$ for
$n \geq dom (s')$).
\smallskip
{\capit Claim.} {\it The following set $A$ is large: $s \in A 
\Longleftrightarrow$ for some $k < \omega$ and $\langle t_\ell , f_\ell^1,
f_\ell^2 ; \; \ell \in \omega \rangle$
we have $s \subseteq t_\ell$, $t_\ell \in \omega^k$, $t_\ell (lh(s))
\geq \ell$, $f_\ell^1 \neq f_\ell^2 \in 2^\omega$, $f_\ell^1 \restrict
\ell = f_\ell^2 \restrict \ell$, and $\forall f \in \omega^\omega$ (with
$t_\ell \subseteq f$)
$\forall m \in \omega \; \forall i \in \{ 1,2 \} \; ( (t_\ell ,f)
\not{\forces} f_\ell^i \restrict m \not\in \breve T )$.}
\smallskip
{\it Proof.} Let $sp \;\breve T$ be the $\DD$-name for the
subset of $\omega$ which describes the levels at which there is
a splitting node in $\breve T$. By thinning out $T$ (in the
generic extension) if necessary, we may assume that
\smallskip
\centerline{$\forces_\DD$ the $j$-th member of $sp \;\breve T$
(denoted by $\tau_j$) is $> \breve d (j)$,}
\smallskip
\noindent where $\breve d$ is (as always) the $\DD$-name for the
Hechler real. Let $(s^* , f^* ) \in \DD$, $lh (s^*) = j^*$. So
$(s^* , f^*)$ forces no bound on $\tau_{j^*}$ --- even no
$(s^* , f')$ does (*). We assume there is no $s \in A$ with
$(s, f^*) \leq (s^* , f^*)$ and reach a contradiction.
\par
Let $I$ be the dense set of conditions forcing a value to
$\tau_{j^*}$; and let $B = \{ s \in \omega^{< \omega} ; \; \exists
f \in \omega^\omega \; ( (s,f) \in I )\}$. By the main lemma 1.2.
we have $rk(s^*, B) < \omega_1$. We prove by induction on the 
ordinal $\beta < \omega_1$
\smallskip
\item{(**)} {\it if $s \in \omega^{< \omega}$ is such that $(s, f^*) \leq
(s^* , f^*)$ and $rk(s,B) = \beta$, then $\exists m < \omega \;
\forall f \in \omega^\omega $ (with $s \subseteq f$) 
$( (s,f) \not{\forces} \tau_{j^*}
\neq m)$.}
\smallskip
\noindent If we succeed for $s = s^*$ then we get a contradiction
to (*).
\par
$\beta = 0$. So $s \in B$. Thus for some $f' \geq f^*$, $(s,f')$ forces
a value to $\tau_{j^*}$: $(s,f') \forces \tau_{j^*} = m$, for some
$m \in \omega$, giving (**).
\par
$\beta > 0$. By the definition of rank there are $k \in \omega$,
$t_\ell \in \omega^k$ ($\ell \in \omega$) such that $s \subseteq t_\ell$,
$t_\ell (lh (s)) \geq \ell$, and $rk (t_\ell , B) = \beta_\ell < \beta$.
(We consider only $\ell$ with $\ell \geq max(rng(f^* \restrict k))$.)
By induction hypothesis there are $m_\ell \in \omega$ such that
$\forall f \in \omega^\omega$ (with $t_\ell \subseteq f$) 
$( (t_\ell , f) \not{\forces} \tau_{j^*}
\neq m_\ell )$. We consider two subcases.
\par
{\sanse Case 1.} For some $m$ we have infinitely many $\ell$
such that $m_\ell = m$. Then we can use this $m$ for $s$ and
get (**).
\par
{\sanse Case 2.} $\langle m_\ell ; \; \ell \in \omega \rangle$
converges to $\infty$. Replacing it by a subsequence, if necessary,
we may assume that it is strictly increasing. We show that 
$\langle t_\ell ; \; \ell \in \omega \rangle$ witnesses $s \in A$,
contradicting our initial assumption. \par
For each $\ell$ let $T_\ell = \{ \rho \in 2^{< \omega} ; \;$
for no $f \in \omega^\omega$ does $(t_\ell , f) \forces \rho \not\in 
\breve T \}$. Clearly $T_\ell \subseteq 2^{< \omega}$, $\langle
\rangle \in T_\ell$, and $T_\ell$ is closed under initial segments.
Also we have that $\rho \in T_\ell$ implies either $\rho \hat{\;}
\langle 0 \rangle \in T_\ell$ or $\rho \hat{\;} \langle 1 \rangle \in
T_\ell$ (otherwise we can find $f_0, f_1 \in \omega^\omega$ such
that $(t_\ell , f_0) \forces \rho \hat{\;} \langle 0 \rangle \not\in
\breve T$ and $(t_\ell, f_1 ) \forces \rho \hat{\;} \langle 1 \rangle
\not\in \breve T$; let $f = max \{ f_0 , f_1 \}$; choose $p \leq
(t_\ell , f)$ such that $p \forces \rho \in \breve T$ (by assumption
on $\rho$); but then there exists $q \leq p$ such that either
$q \forces \rho \hat{\;} \langle 0 \rangle \in \breve T$ or
$q \forces \rho \hat{\;} \langle 1 \rangle \in \breve T$, a contradiction).
\par
Finally, $T_\ell$ has a splitting node at level $m_\ell$; i.e. for
some $\rho = \rho_\ell \in T_\ell \cap 2^{m_\ell}$, we have $\rho
\hat{\;} \langle 0 \rangle \in T_\ell$ and $\rho \hat{\;} \langle 1 
\rangle \in T_\ell$ (if not, for each $\rho \in 2^{m_\ell} \;
\exists f_\rho \in \omega^\omega$ such that $(t_\ell , f_\rho)
\forces "\rho \hat{\;} \langle 0 \rangle \not\in \breve T$ or $\rho
\hat{\;} \langle 1 \rangle \not\in \breve T"$; let $f = max
\{ f_\rho ; \; \rho \in 2^{m_\ell} \}$. We know that $(t_\ell , f)
\not{\forces} m_\ell \neq \tau_{j^*}$; so there is $p \leq (t_\ell,
f)$ such that $p \forces m_\ell = \tau_{j^*}$; i.e. $p \forces
m_\ell \in sp \;\breve T$; we now get a contradiction as before).
\par
Hence we can find $f_\ell^1, f_\ell^2 \in [T_\ell]$ such that
$f^1_\ell \restrict (m_\ell + 1) = \rho_\ell \hat{\;} \langle
0 \rangle$ and $f_\ell^2 \restrict (m_\ell + 1) = \rho_\ell \hat{\;}
\langle 1 \rangle$. Thus $\langle t_\ell ; \; \ell \in \omega \rangle$,
$\langle f_\ell^1, f_\ell^2 ; \; \ell \in \omega \rangle$ witness
$s \in A$. This final contradiction proves the claim. $\qed$
\smallskip
{\it Continuation of the proof of the theorem.}
We assume that $\forces_\DD \breve T = \{ \tau_j ; \; j \in \omega \}$;
i.e. $\tau_j [G]$ ($j \in \omega$) will enumerate the tree $T = \breve
T [G]$ in the generic extension. We also let $\breve T_j$ be
the name for the tree $T_{\tau_j [G]}$; i.e. $\forces_\DD
\breve T_j = \{ \nu \in \breve T ; \; \nu \subseteq \tau_j$ or $\tau_j
\subseteq \nu \}$. For each $j \in \omega$ there is --- according
to the claim for $\breve T_j$ instead of $\breve T$ --- 
a large set $A_j \subseteq \omega^{<\omega}$; and
for $s \in A_j$ there is a sequence $\langle t_\ell^{s,j} , f_\ell^{
1,s,j} , f_\ell^{2,s,j} ; \; \ell \in \omega \rangle$ that
witnesses $s \in A_j$. For every $j \in \omega$, $s \in A_j$
and $m \in \omega$ we define $S_{j,s,m} = \{ f_\ell^{i,s,j} \restrict
k ; \; k \in \omega, \; i \in \{1,2\} , \; m \leq \ell\in \omega
\}$. By construction the function $f_{j,s,m}$ defined by
$f_{j,s,m} (k) = \vert S_{j,s,m} \cap 2^k \vert$ converges to
$\infty$. By assumption 1) we can choose $g \in \omega^\omega 
\cap M$ such that $\forall j , s, m \; \exists^\infty i \;
(\vert S_{j,s,m} \cap 2^{g(i)} \vert > 4^i$). \par
Now let $U$ be as in the proof of theorem 3.1.; and choose $u^* \in
U$ as there (i.e. $u^*$ is random over a countable model $N$
containing $g$ and all $S_{j,s,m}$ --- using assumption 2)).
We also define $B$ and $B_k$ ($k \in \omega$) as in the proof of theorem 
3.1.
\par
We assume that $\forces_\DD " \breve T$ is a perfect set of reals
random over $V"$; in particular $\forces_\DD \breve T \subseteq
\bigcup_{k \in \omega} B_k$. So there are $(s^* , f^*) \in \DD$,
$j \in \omega$ and $k \in \omega$ such that
\smallskip
\centerline{$(s^*, f^*) \forces_\DD \breve T_j \subseteq B_k$}
\smallskip
\noindent (cf the corresponding argument in the proof of theorem 3.1.).
Without loss $s^* \in A_j$ (otherwise increase the condition using
the claim). Let $m > max (rng (f^* \restrict k_{j, s^*}))$ where
$k_{j, s^*}$ is such that for all $\ell \in \omega$, 
$t_\ell^{s^* ,j} \in\omega^{k_{j,s^*}}$.
Then $\forall \ell \geq m$, $(t_\ell^{s^*,j}, f^*)$ is an extension
of $(s^* , f^*)$. So we must have $S_{j,s^*,m} \subseteq B_k$
(because for any element of the former set we have an extension
of $(s^*, f^*)$ forcing this element into $\breve T_j$).
\par
The rest of the proof is again as in the proof of theorem 3.1.
For infinitely many $i$ we have $\vert S_{j,s^*,m} \cap 2^{g(i)}
\vert > 4^i$; for each such $i$, the probability that $u^* (i)
\cap ( S_{j, s^* , m} \cap 2^{g(i)}) = \emptyset$ is $\leq
e^{-2^i}$; the probability that this happens infinitely often
is zero, contradicting the fact that $u^*$ is random over $N$.
$\qed$
\Bigskip

\noindent{\dunhg References}
\Smallskip
\itemitem{[Ba]} {\capit T. Bartoszy\'nski,} {\it Combinatorial
aspects of measure and category,} Fundamenta Mathematicae, vol. 127
(1987), pp. 225-239.
\smallskip
\itemitem{[BaJ]} {\capit T. Bartoszy\'nski and H. Judah,} {\it
Jumping with random reals,} Annals of Pure and Applied Logic,
vol. 48 (1990), pp. 197-213.
\smallskip
\itemitem{[BJS]} {\capit T. Bartoszy\'nski, H. Judah and S.
Shelah,} {\it The Cicho\'n diagram,} submitted to Journal of Symbolic
Logic.
\smallskip
\itemitem{[BD]} {\capit J. Baumgartner and P. Dordal,} {\it
Adjoining dominating functions,} {Journal of Symbolic Logic,} vol.
50 (1985), pp. 94-101.
\smallskip
\itemitem{[BrJ]} {\capit J. Brendle and H. Judah,} {\it Perfect
sets of random reals,} submitted to Israel Journal of Mathematics.
\smallskip
\itemitem{[CP]} {\capit J. Cicho\'n and J. Pawlikowski,} {\it
On ideals of subsets of the plane and on Cohen reals,}
Journal of Symbolic Logic, vol. 51 (1986), pp. 560-569.
\smallskip
\itemitem{[Fr]} {\capit D. Fremlin,} {\it Cicho\'n's diagram,}
S\'eminaire Initiation \`a l'Analyse (G. Choquet,
M. Rogalski, J. Saint Raymond), Publications Math\'ematiques
de l'Universit\'e Pierre et Marie Curie, Paris, 1984,
pp. 5-01 - 5-13. 
\smallskip
\itemitem{[GS]} {\capit M. Gitik and S. Shelah,} {\it More on ideals with 
simple forcing notions,} to appear in Annals of Pure and Applied Logic.
\smallskip
\itemitem{[Je]} {\capit T. Jech,} {\it Set theory,} Academic Press,
San Diego, 1978.
\smallskip
\itemitem{[Ju]} {\capit H. Judah,} {\it Absoluteness for projective
sets,} to appear in Logic Colloquium 1990.
\smallskip
\itemitem{[JS]} {\capit H. Judah and S. Shelah,} {\it The Kunen-Miller
chart (Lebesgue measure, the Baire property, Laver reals and
preservation theorems for forcing),} Journal of
Symbolic Logic, vol. 55 (1990), pp. 909-927.
\smallskip
\itemitem{[Ku]} {\capit K. Kunen,} {\it Set theory,} North-Holland,
Amsterdam, 1980.
\smallskip
\itemitem{[Mi]} {\capit A. Miller,} {\it Arnie Miller's problem list,}
to appear in Proccedings of the Bar-Ilan conference on set theory
of the reals, 1991.
\smallskip
\itemitem{[Pa]} {\capit J. Pawlikowski,} {\it Why Solovay
real produces Cohen real,} Journal of Symbolic Logic,
vol. 51 (1986), pp. 957-968.
\smallskip
\itemitem{[Tr]} {\capit J. Truss,} {\it Sets having calibre 
$\aleph_1$,} Logic Colloquium 76, North-Holland, Amsterdam,
1977, pp. 595-612.

\vfill\eject\end